\documentclass[11pt, letterpaper]{article}

\pdfoutput=1

\usepackage{amsmath}
\usepackage{amsfonts}
\usepackage{amssymb} 
\usepackage{amsthm}
\usepackage{mathrsfs} 
\usepackage{mathtools} 
\usepackage[nottoc,numbib]{tocbibind} 
\usepackage[usenames,dvipsnames,svgnames,table]{xcolor} 

\usepackage{tikz}
\usetikzlibrary{shapes,matrix,arrows,chains,trees,calc} 

\usepackage{xr} 
\usepackage{cite} 
\usepackage{enumerate} 
\usepackage[colorlinks=true, allcolors=MidnightBlue, pdftitle={5-dimensional geometries I: the general classification}, pdfauthor={Andrew Geng}]{hyperref}
\usepackage{url}
\usepackage[total={6.5in,9in}]{geometry}
\usepackage[T1]{fontenc} 
\usepackage[utf8]{inputenc} 

\title{5-dimensional geometries I: the general classification}
\author{Andrew Geng}

\theoremstyle{plain}
\newtheorem{thm}{Theorem}[section]

\theoremstyle{definition}
\newtheorem{defn}[thm]{Definition}

\theoremstyle{remark}

\makeatletter
\let\c@figure\c@thm
\let\c@table\c@thm
\makeatother
\numberwithin{figure}{section}
\numberwithin{table}{section}

\newcommand{\keyword}{\emph}
\newcommand{\op}{\operatorname}

\newcommand{\R}{\mathbb{R}}
\newcommand{\C}{\mathbb{C}}
\newcommand{\Q}{\mathbb{Q}}

\newcommand{\Z}{\mathbb{Z}}
\newcommand{\Hyp}{\mathbb{H}}
\newcommand{\Euc}{\mathbb{E}}
\newcommand{\Heis}{\mathrm{Heis}}
\newcommand{\Sol}{\mathrm{Sol}}
\newcommand{\Nil}{\mathrm{Nil}}
\newcommand{\SLcover}{\widetilde{\mathrm{SL}_2}}

\newcommand{\lie}{\mathfrak}
\newcommand{\tanalg}{T_{\mathbf{1}}} 

\newcommand{\SemiR}[2]{ { \underset{ {#2} }{ { {#1} } \rtimes \R} } }

\begin{document}

\maketitle

\begin{abstract}
This paper is the first of a 3-part series that classifies the
$5$-dimensional Thurston geometries.
The present paper (part 1 of 3) summarizes the general classification,
giving the full list, an outline of the method, and some illustrative examples.
This includes phenomena that have not appeared in lower dimensional geometries,
such as an uncountable family of geometries $\SLcover \times_\alpha S^3$.
\end{abstract}

\tableofcontents

\section{Introduction}

By the classification of closed surfaces
(see e.g. \cite[Thm.\ 5.11]{matsumoto_morse}),
every closed surface is diffeomorphic to a quotient
of $\Euc^2$, $S^2$, or $\Hyp^2$ by a discrete group of isometries.
It is a classical result that in dimension $2$,
these three are the only connected, simply-connected, complete Riemannian manifolds
with transitive isometry group
(see e.g.\ \cite[Thm.\ 3.8.2]{thurstonbook}).

The quest for the $3$-dimensional generalization
that became Thurston's Geometrization Conjecture led to
a version of the following definition.
(The equivalence to older definitions is outlined in
Part II, \cite[Prop.~\ref{ii:prop:geometries3}]{geng2}.)
\begin{defn}[\textbf{Geometries}, following {\cite[Defn.\ 3.8.1]{thurstonbook}} and {\cite[\S 1.1]{filipk}}] \label{defn:geometries1} ~
    \begin{enumerate}[(i)]
        \item A \keyword{geometry} is a connected, simply-connected homogeneous
            space $M = G/G_p$ where $G$ is a connected Lie group acting faithfully
            with compact point stabilizers $G_p$.
        \item $M$ is a \keyword{model geometry} if there is some lattice
            $\Gamma \subset G$ that acts freely on $M$.
            Then the manifold $\Gamma \backslash G / G_p$
            is said to be \keyword{modeled on} $M$.
        \item $M$ is \keyword{maximal} if it is not $G$-equivariantly
            diffeomorphic to any other geometry $G'/G'_p$ with $G \subsetneq G'$.
            Any such $G'/G'_p$ is said to \keyword{subsume} $G/G_p$.
    \end{enumerate}
\end{defn}
Then a closed $3$-manifold is a quotient
of at most one maximal model geometry, which can be determined
from the fundamental group \cite[Thm.\ 4.7.8]{thurstonbook} or
from the existence of certain bundle structures (usually Seifert bundles)
and some topological data (usually two Euler numbers) \cite[Thm.\ 5.3]{scott}.
Thurston classified the $3$-dimensional maximal model geometries
and found eight (see \cite[Thm.\ 3.8.4]{thurstonbook}).

In $4$ dimensions, Filipkiewicz classified the maximal model geometries in \cite{filipk}.
Though $4$-manifolds without geometric decompositions \cite[\S{13.3} \#3]{hillman}
indicate there is less hope for a straightforward generalization of geometrization,
Filipkiewicz's classification highlights a few interesting firsts.
The list comprises $18$ geometries and---for
the first time---a countably infinite family,
named $\Sol^4_{m,n}$. (See e.g.\ \cite[\S 7.1]{hillman} or
\cite[\S 1, Table 1]{wall86} for the names currently in use.)
One of the eighteen is $\mathbb{F}^4 = \R^2 \rtimes \op{SL}(2,\R)/\op{SO}(2)$,
the first geometry
to admit finite-volume quotients but no compact quotients.


The direction to take should now seem straightforward.
One seeks a classification of maximal model geometries in all dimensions;
but a handful of obstacles stand in the way of such a classification:
\begin{enumerate}
    \item Existing classifcations, including now the present paper,
        rely on tools that may become unusable with increasing dimension.
        For example, the case of discrete point stabilizers
        (\hspace{1sp}\cite[Thm.\ 3.8.4(c)]{thurstonbook} in dimension $3$,
        \cite[Ch.\ 6]{filipk} in dimension $4$,
        and \cite[Thm.~\ref{ii:thm:main}(ii)]{geng2} in Part II)
        relies on a classification of solvable Lie algebras over $\R$,
        which is incomplete in dimensions $7$ and up.
        (See e.g.\ \cite[Introduction]{snobl2012} for a summary of known progress,
        and \cite{boza} for a wider survey.)
    \item The aforementioned aspects of the $4$-dimensional classification
        suggest that new phenomena may continue to appear for a few more dimensions.
        A workable approach to a general classification may not be evident without
        knowledge of such features.
\end{enumerate}
An optimistic interpretation of these obstacles is that the $5$-dimensional case
is both tractable and potentially illustrative.
Having carried out the classification, the new phenomena
are summarized in Section \ref{sec:whatsnew};
the main result is the following list.

\begin{thm}[\textbf{Classification of 5-dimensional geometries}] \label{thm:full_list}
    The maximal model geometries of dimension $5$ are:
    \begin{enumerate}
        \item The geometries with constant curvature:
            \begin{align*}
                \Euc^5 &= \R^5 \rtimes \op{SO}(5) / \op{SO}(5) &
                S^5 &= \op{SO}(6)/\op{SO}(5) &
                \Hyp^5 &= \op{SO}(5,1)/\op{SO}(5) ;
            \end{align*}
        \item The other irreducible Riemannian symmetric spaces
            $\op{SL}(3,\R)/\op{SO}(3)$ and $\op{SU}(3)/\op{SO}(3)$;
        \item The unit tangent bundles or universal covers of circle bundles:
            \begin{align*}
                T^1(\Hyp^3) &= \op{PSL}(2,\C)/\op{SO}(2)  &
                T^1(\Euc^{1,2}) &= \R^3 \rtimes \op{SO}(1,2)^0 / \op{SO}(2)  &
                & \widetilde{\op{U}(2,1)/\op{U}(2)} ;
            \end{align*}
        \item The associated bundles 
            (see e.g.\ \cite[\S{1.3} Vector Bundles]{sharpe} for the notation):
            \begin{align*}
                \Heis_3 \times_\R S^3
                    &= (\Heis_3 \rtimes \widetilde{\op{SO}(2)}) \times S^3 /
                        \{(0,0,s), \gamma(t), e^{\pi i s} \}_{s,t \in \R} \\
                \Heis_3 \times_\R \SLcover
                    &= (\Heis_3 \rtimes \widetilde{\op{SO}(2)}) \times \SLcover /
                        \{(0,0,s), \gamma(t), \gamma(s) \}_{s,t \in \R} \\
                \SLcover \times_\alpha S^3
                    &= \SLcover \times S^3 \times \R /
                        \{\gamma(s), e^{\pi i t}, \alpha s + t \}_{s, t \in \R},
                        \quad 0 < \alpha < \infty \\
                \SLcover \times_\alpha \SLcover
                    &= \SLcover \times \SLcover \times \R /
                        \{\gamma(s), \gamma(t), \alpha s + t \}_{s, t \in \R},
                        \quad 0 < \alpha \leq 1 \\
                L(a;1) \times_{S^1} L(b;1)
                    &= S^3 \times S^3 \times \R /
                        \{e^{\pi i s}, e^{\pi i t}, as + bt\}_{s, t \in \R},
                        \quad 0 < a \leq b \text{ coprime in } \Z ,
            \end{align*}
            where the Heisenberg group $\Heis_3$ is $\R^3$ with the multiplication law
                \[ (x,y,z)(x',y',z') = (x+x', y+y', z+z'+xy'-x'y), \]
            on which $\op{SO}(2)$ acts through the action of $\op{SL}(2,\R)$
            on the $x,y$ plane, and $t \mapsto e^{\pi i t} \in S^3$
            and $\gamma: \R \to \widetilde{\op{SO}(2)} \subset \SLcover$
            are $1$-parameter subgroups sending $\Z$ to the center;
        \item The three principal $\R$-bundles with non-flat connections
            over the $\mathbb{F}^4$ geometry,
            distinguished from each other by their curvatures:
            \begin{align*}
                \R^2 \rtimes \SLcover
                    &\cong (\R^2 \rtimes \SLcover) \rtimes \op{SO}(2) / \op{SO}(2) \\
                \mathbb{F}^5_a
                    &= \Heis_3 \rtimes \SLcover / \{ (0,0,at), \gamma(t) \}_{t \in \R} ,
                    \quad a = 0 \text{ or } 1 ;
            \end{align*}
        \item The six simply-connected indecomposable nilpotent Lie groups,
            named by their Lie algebras as in \cite[Table II]{patera},
            in which the point stabilizer of the identity element
            is a maximal compact group of automorphisms
            (specified in Table \ref{table:geoms_by_isotropy}):
            \begin{align*}
                A_{5,1} &= \SemiR{\R^4}{x^2,\, x^2}  &
                A_{5,2} &= \SemiR{\R^4}{x^4}  &
                A_{5,3} &= \SemiR{(\R \times \Heis_3)}{x_3 \to x_2 \to y}  \\
                A_{5,4} &= \Heis_5  &
                A_{5,5} &= \SemiR{\Nil^4}{3 \to 1}  &
                A_{5,6} &= \SemiR{\Nil^4}{4 \to 3 \to 1} ;
            \end{align*}
        \item The simply-connected indecomposable non-nilpotent solvable Lie groups,
            specified the same way:
            \begin{align*}
                A^{a,b,-1-a-b}_{5,7} &= \SemiR{\R^4}{\text{4 distinct real roots}}  &
                A^{1,-1-a,-1+a}_{5,7} &= \SemiR{\R^4}{\text{2 complex, 2 distinct real}}  \\
                A^{1,-1,-1}_{5,7} &= \SemiR{\R^4}{x-1,\, x-1,\, x+1,\, x+1}  &
                A^{-1}_{5,8} &= \SemiR{\R^4}{x^2,\, x-1,\, x+1}  \\
                A^{-1,-1}_{5,9} &= \SemiR{\R^4}{(x-1)^2,\, x+1,\, x+1}  &
                A^{-1}_{5,15} &= \SemiR{\R^4}{(x-1)^2,\, (x+1)^2}  \\
                A^{0}_{5,20} &= \SemiR{(\R \times \Heis_3)}{\text{Lorentz},\, y \to x_1} &
                A^{-1,-1}_{5,33} &= \R^3 \rtimes \{xyz=1\}^0 ;
            \end{align*}
        \item and all twenty-nine products of lower-dimensional geometries
            involving no more than one Euclidean factor,
            named as in \cite[Table 1]{wall86}.
            \begin{enumerate}
                \item $4$-by-$1$:
                    \begin{align*}
                        S^4 &\times \Euc  &  \Hyp^4 &\times \Euc  &
                        \mathbb{CP}^2 &\times \Euc  &  \mathbb{CH}^2 &\times \Euc  &
                        \mathbb{F}^4 &\times \Euc  \\
                        \Nil^4 &\times \Euc  &
                        \Sol^4_0 &\times \Euc  &  
                        \Sol^4_1 &\times \Euc  &
                        \Sol^4_{m,n} &\times \Euc
                    \end{align*}
                \item $3$-by-$2$:
                    \begin{align*}
                        & & \Euc^3 &\times S^2 &
                        \Euc^3 &\times \Hyp^2 \\
                        S^3 &\times \Euc^2 &
                        S^3 &\times S^2 &
                        S^3 &\times \Hyp^2 \\
                        \Hyp^3 &\times \Euc^2 &
                        \Hyp^3 &\times S^2 &
                        \Hyp^3 &\times \Hyp^2 \\
                        \Heis_3 &\times \Euc^2 &
                        \Heis_3 &\times S^2 &
                        \Heis_3 &\times \Hyp^2 \\
                        \Sol^3 &\times \Euc^2 &
                        \Sol^3 &\times S^2 &
                        \Sol^3 &\times \Hyp^2 \\
                        \SLcover &\times \Euc^2 &
                        \SLcover &\times S^2 &
                        \SLcover &\times \Hyp^2
                    \end{align*}
                \item $2$-by-$2$-by-$1$:
                    \begin{align*}
                        S^2 &\times S^2 \times \Euc &
                        S^2 &\times \Hyp^2 \times \Euc &
                        \Hyp^2 &\times \Hyp^2 \times \Euc
                    \end{align*}
            \end{enumerate}
    \end{enumerate}
\end{thm}
More explicit instructions for constructing these geometries---such as
the solvable Lie groups and their automorphism groups---are
delegated to where they occur in the classification in Parts II and III
\cite{geng2,geng3}.

\paragraph{Roadmap.}
The present paper (Part I) summarizes the classification;
Section \ref{sec:whatsnew} picks out illustrative examples,
Section \ref{sec:method} outlines the strategy,
and Section \ref{sec:whatelse} briefly surveys related classifications.
Part II \cite{geng2} classifies the point stabilizer subgroups $G_p$
and classifies the geometries where $G_p$ acts irreducibly or trivially
on tangent spaces.
Part III \cite{geng2} classifies the remaining geometries
after showing that they all admit invariant fiber bundle structures
(hence the name ``fibering geometries'').


\paragraph{Acknowledgments.}

I have the pleasure to thank my advisor Benson Farb
for years of helpful discussions, patient advice,
and extensive comments on drafts of all three papers in this series.
Thanks also to my second advisor Danny Calegari,
especially for his keen sense of where the interesting unanswered
questions lay;
to Jonathan Hillman, Christoforos Neofytidis,
Daniel Studenmund, Ilya Grigoriev, and Nick Salter for
a number of other enlightening conversations;
to Ilka Agricola for pointing out related work on $\op{SO}(3)_5$-structures;
and to Greg Friedman, Christoforos Neofytidis, and Kenneth Knox
for opportunities to speak about an early version of these results.
Finally, I would like to thank the University of Chicago for
support during this work.

\section{Salient examples} \label{sec:whatsnew}

\subsection{New phenomena}

Much of our interest in this classification is in the search for phenomena
that occur for the first time in dimension $5$,
in hopes of finding a pattern that continues in higher dimensions.
See also Section \ref{sec:method} for a discussion of new tools.

\paragraph{An uncountable family of geometries.}
    The associated bundles
    $\SLcover \times_\alpha S^3$ ($0 < \alpha < \infty$) form
    an \emph{uncountable} family of maximal model geometries.
    (Taking $\SLcover \times_\alpha S^3$ as a circle bundle over $S^2 \times \Hyp^2$
    with an invariant connnection,
    the parameter $\alpha$ is a ratio of curvatures
    in the $S^2$ and $\Hyp^2$ directions.)
    This and $\SLcover \times_\alpha \SLcover$
    are the first occurrences of uncountable families.
    Since every lattice in a Lie group is finitely presented
    \cite[Thm.~I.1.3.1]{onishchik2},
    $\pi_1$ of the quotient manifolds will not determine the geometries.
    Details are in Part III,
    \cite[Prop.~\ref{iii:prop:fiber2_associated_lattice}]{geng3}.

\paragraph{An infinite family without compact quotients.}
    In fact, $\SLcover \times_\alpha S^3$
    admits compact quotients if and only if $\alpha$ is rational
    \cite[Prop.~\ref{iii:prop:fiber2_compact_rationals}]{geng3}.
    (Recall that beginning with $\Hyp^2$, geometries can have noncompact quotients
    of finite volume; and beginning with $\mathbb{F}^4 = \R^2 \rtimes \op{SL}(2,\R)/\op{SO}(2)$,
    model geometries might have no compact quotients.)

\paragraph{Non-unique maximality.}
    The geometry $T^1 S^3 = \op{SO}(4)/\op{SO}(2)$ is a non-maximal form of
    $S^3 \times S^2$ and $L(1;1) \times_{S^1} L(1;1)$---both of which are
    maximal \cite[Rmk.~\ref{iii:rmk:nonunique_maximality}]{geng3}.
    This contrasts with the positive results for unique maximality
    listed in the discussion after \cite[Prop.~1.1.2]{filipk}.

\paragraph{Inequivalent compact geometries with the same diffeomorphism type.}
    Using Barden's diffeomorphism classification \cite{barden}
    of simply-connected $5$-manifolds by second homology
    and second Stiefel-Whitney class, one can prove that
    the associated bundles of lens spaces $L(a;1) \times_{S^1} L(b;1)$
    are all diffeomorphic to $S^3 \times S^2$ \cite[Cor.~3.3.2]{ottenburger2009}.

    More broadly one can attempt to give the classification up to diffeomorphism,
    following previous results such as \cite[Cor.~p.~624]{mostow1950},
    \cite{gorbatsevich1977}, \cite{ishihara1955}, and \cite[Thm.~1.0.3]{ottenburger2009}.
    Most of the geometries are products of $\R^k$ and some spheres;
    the two exceptions are
    $\mathbb{CP}^2 \times \Euc$ and the rational homology sphere
    $\op{SU}(3)/\op{SO}(3)$, named $X_{-1}$ in Barden's classification
    \cite[Introduction]{boyergalicki}.

    Note that while the correct diffeomorphism type may be obvious enough to guess,
    it is not as obviously correct.
    The problem is proving that the space is a direct product of $\R^k$
    and the product of spheres onto which it deformation retracts---such
    a claim is false for any nontrivial vector bundle (such as $TS^2$)
    and for a homogeneous example by Samelson
    discussed in \cite[\S{5} Example 4]{mostow1955}.
    Instead one has to use either an explicit description of the diffeomorphism type
    from \cite[Thm.~A]{mostow_covariant_1962}
    or the fact that sufficiently nice bundles over contractible spaces are trivial
    \cite[Cor.~10.3]{husemoller}.

\paragraph{Isotropy irreducible spaces.}
    The geometries $\op{SL}(3,\R)/\op{SO}(3)$ and $\op{SU}(3)/\op{SO}(3)$
    have point stabilizers $\op{SO}(3)$ acting irreducibly on the ($5$-dimensional)
    tangent spaces.
    This is the first occurrence of such an action 
    being irreducible and not the standard representation
    for a group of the same isomorphism type.
    These two geometries are still symmetric spaces, but in higher dimensions there exist
    homogeneous spaces with irreducibly-acting point stabilizers that are not
    symmetric spaces (See e.g.\ \cite[Introduction]{heintze1996}).

\subsection{Examples that highlight tools}

The classification of geometries requires an increasingly wide range of tools
as the dimension increases. These are a handful of examples where
either unexpected tools appeared, or familiar tools exhibit behavior
that is not completely obvious at first glance.

\paragraph{Model geometries via Galois theory and Dirichlet's unit theorem.}
    Some Galois theory is needed to answer questions of lattice existence,
    such as to prove that $\R \rtimes \op{Conf}^+ \Euc^3 / \op{SO}(3)$
    (where the action of $\op{Conf}^+ \Euc^3$ on $\R$ is chosen to make
    the semidirect product unimodular) is not a model geometry
    \cite[Prop.~\ref{iii:prop:main:ii}(iv)]{geng3}.

    Dirichlet's unit theorem makes an appearance when we construct
    a lattice in $\R^3 \rtimes \{xyz = 1\}^0$ by taking a finite index subgroup of
    $\mathcal{O}_K \rtimes \mathcal{O}_K^\times$ where $K$
    is a totally real cubic field extension of $\Q$
    \cite[Prop.~\ref{ii:prop:commutingmatrices}]{geng2}.

\paragraph{Point stabilizers not realized.}
    The classification of geometries starts by classifying subgroups of $\op{SO}(5)$
    in order to classify point stabilizers---but not every subgroup is realized by
    a maximal model geometry. For example, $\op{SO}(3)$ in its standard representation
    is one such subgroup, though the non-model geometry
    $\R \rtimes \op{Conf}^+ \Euc^3/\op{SO}(3)$
    mentioned above suggests this can be thought of as a near miss.
    The non-occurrence of $\op{SU}(2)$ is another example,
    a feature shared by the $4$-dimensional classification of geometries.
    Other subgroups---namely $\op{SO}(4)$ and $\op{SO}(3) \times \op{SO}(2)$---are
    point stabilizers only of product geometries.
    A listing of (non-product) geometries by point stabilizer
    is given below in Table \ref{table:geoms_by_isotropy}
    after Figure \ref{fig:isotropy_poset} names the subgroups.

\paragraph{Geometries in higher dimensions with reducible isotropy and no fibering.}
    When point stabilizers act reducibly on tangent spaces,
    our strategy breaks down the problem by showing the existence
    of an invariant fiber bundle structure.
    That this is possible is is a convenient accident of low dimensions;
    higher dimensions introduce
    isotropy-reducible geometries that admit no fibering.

    For example, in dimension $18$, there is $\op{Sp}(3)/\op{Sp}(1)$,
    where the embedding $\op{Sp}(1) \hookrightarrow \op{Sp}(3)$
    is given by the irreducible representation of
    $\op{Sp}(1) \cong \op{SU}(2)$ on $\C^6$.
    This has two isotropy summands but admits no nontrivial fibering
    since $\op{Sp}(1)$ is maximal
    (so no larger group can be a point stabilizer of the base space)
    \cite[Example V.10]{dickinsonkerr}.
    A strategy that continues to break the problem down using invariant
    fiber bundle structures may have to account for these exceptions,
    likely using Dynkin's work on classifying maximal subgroups
    of semisimple Lie groups in \cite{dynkin2000maximal, dynkin2000semisimple}.

\paragraph{Non-geometries as base spaces of fiber bundles.}
    Even when invariant fiber bundle structures exist,
    a number of compliactions prevent the classification from
    having a straightforward recursive solution.
    Filipkiewicz warns in \cite[Prop.~2.1.3]{filipk}
    that the base space of an invariant fiber bundle structure
    may fail to be a geometry due to noncompact point stabilizers---e.g.\ the
    action of $\op{PSL}(2,\C)$ on $S^2 \cong \mathbb{CP}^1$
    makes $T^1 \Hyp^3 = \op{PSL}(2,\C)/\op{PSO}(2)$ a fiber bundle over $S^2$.
    Even when point stabilizers are compact, the base may fail to
    be maximal (e.g.\ $\Heis_5$ fibers over $\Euc^4$ with
    $\op{U}(2)$ point stabilizers)
    or a model geometry (e.g.\ $\op{Sol}^3$ over $\op{Aff}^+ \R$,
    which cannot admit a lattice since it is not unimodular).

\subsection{Examples that clarify how the classification is organized}

The eight categories of Theorem \ref{thm:full_list} and the grouping of
geometries into parametrized families involved some
arbitrary choices.
This section discusses the chosen method of organization and some variations.

\paragraph{The omission of some spaces that one might have guessed.}
    Some of the categories in Thm.~\ref{thm:full_list}
    are conspicuously missing geometries
    that happen to be non-model or non-maximal.
    \begin{enumerate}
        \item[3.]
            The tautological unit circle bundle $\op{U}(3)/\op{U}(2)$
            over $\mathbb{CP}^2$ is non-maximal, being equivariantly
            diffeomorphic to $S^5$.
            The two other unit tangent bundles of 3-dimensional
            spaces of constant curvature are also non-maximal.
            \begin{align*}
                T^1(S^3) = \op{SO}(4)/\op{SO}(2) &\cong S^2 \times S^3  &
                T^1(\Euc^3) = \R^3 \rtimes \op{SO}(3) / \op{SO}(2) &\cong S^2 \times \Euc^3
            \end{align*}
        \item[7.] Many of the solvable Lie groups arising from the list
            in \cite[Table II]{patera} are not unimodular and hence do not
            admit lattices.
        \item[8.] Every product geometry with multiple Euclidean factors
            is non-maximal---but all other products are maximal,
            usually as a consequence of the de Rham decomposition theorem.
            (See \cite[Prop.~\ref{iii:prop:products_are_maximal}]{geng3} in Part III.)
    \end{enumerate}

\paragraph{Counting the geometries and families.}
    The list given in Thm.~\ref{thm:full_list} includes 53 individual geometries
    and the following 6 infinite families of geometries.
    \begin{align*}
        L(a;1) \times_{S^1} L(b;1),
            &\quad a \leq b \text{ coprime positive integers} \\
        \SLcover \times_\alpha S^3,
            &\quad 0 < \alpha < \infty \\
        \SLcover \times_\alpha \SLcover,
            &\quad 0 < \alpha \leq 1 \\
        \SemiR{\R^4}{e^{tA}},
            &\quad e^A \text{ semisimple integer matrix with 4 real eigenvalues} \\
        \SemiR{\R^4}{e^{tA}},
            &\quad e^A \text{ semisimple integer matrix with 2 real eigenvalues} \\
        \Sol^4_{m,n} \times \Euc,
            &\quad m, n \in \Z
    \end{align*}
    To some extent, this count depends on interpretation.
    The first three families could be expanded to include
    products of spheres and hyperbolic spaces,
    while the last three could be unified with
    $\SemiR{\R^4}{x-1,\, x-1,\, x+1,\, x+1}$
    to form one large family with a name like $\Sol^5_{m,n,p}$
    (where $m$, $n$, and $p$ are the middle coefficients of the characteristic
    polynomial of $e^A$).
    Indeed, \cite[Table II]{patera} suggests this latter unification
    by listing all semidirect products $\R^4 \rtimes \R$
    with diagonalizable action under the family $A_{5,7}$.
    We keep the subfamilies of $\Sol^5_{m,n,p}$ separate since
    their point stabilizers have different dimensions.

\section{Overview of method}
\label{sec:method}

\begin{figure}[h!]
    \caption[Closed connected subgroups of $\op{SO}(5)$.
        (Duplicate of \ref{fig_thesis:isotropy_poset})]{
        Closed connected subgroups of $\op{SO}(5)$, with inclusions.
        $\op{SO}(3)_5$ denotes $\op{SO}(3)$ acting on its $5$-dimensional
        irreducible representation;
        and $S^1_{m/n}$ acts as on the direct sum $V_m \oplus V_n \oplus \R$
        where $S^1$ acts irreducibly on $V_m$ with kernel of order $m$.
        See Part II, \cite[Prop.~\ref{ii:isotropy_classification}]{geng2}
        for the proof.
    }
    \label{fig:isotropy_poset}
    \begin{center}\begin{tikzpicture}
        \draw (0, 5) node (so5) {$\op{SO}(5)$};
        \draw (-2, 4) node (so4) {$\op{SO}(4)$};
        \draw (1, 4) node (so3so2) {$\op{SO}(3) \times \op{SO}(2)$};
        \draw (4, 4) node (so35) {$\op{SO}(3)_5$};

        \draw (-3, 3) node (u2) {$\op{U}(2)$};
        \draw (-3, 2) node (su2) {$\op{SU}(2)$};

        \draw (2, 3) node (so3) {$\op{SO}(3)$};
        \draw (0, 2) node (t2) {$\op{SO}(2) \times \op{SO}(2)$};

        \draw (-2, 1) node (s11) {$S^1_1$};
        \draw (0, 1) node (s1q) {$S^1_{m/n}$};
        \draw (2, 1) node (so2) {$S^1_0 = \op{SO}(2)$};
        \draw (4, 1) node (s12) {$S^1_{1/2}$};

        \draw (0, 0) node (triv) {$\{1\}$};

        \draw (so5) -- (so4) -- (u2) -- (su2) -- (s11) -- (triv);
        \draw (so5) -- (so3so2) -- (so3) -- (so2) -- (triv);
        \draw (so5) -- (so35) -- (s12) -- (triv);
        \draw (so4) -- (so3);
        \draw (so3so2) -- (t2) -- (so2);
        \draw (u2) -- (t2) -- (s11);
        \draw (t2) -- (s1q) -- (triv);
        \draw (t2) -- (s12);
    \end{tikzpicture}\end{center}
\end{figure}
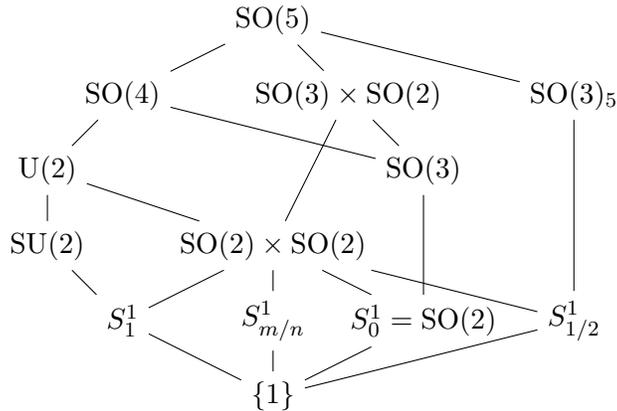

The classification of $5$-dimensional geometries $M = G/G_p$ begins,
following Thurston \cite[\S 3.8]{thurstonbook}
and \cite[\S 1.2]{filipk},
by using the representation theory of compact groups to list the subgroups
$G_p \subseteq \op{SO}(T_p M)$ that could be point stabilizers
(Figure \ref{fig:isotropy_poset}).

\begin{table}[h!]
    \caption[Non-product geometries listed by point stabilizer.]{Using the classification,
        non-product geometries can be listed by point stabilizer.}
    \label{table:geoms_by_isotropy}
    \begin{center}\begin{tabular}{c@{\hskip 24pt}l}
        Stabilizer & Geometries \\
        \hline
        \rule[-9pt]{0pt}{27pt}
        $\op{SO}(5)$ &
        $\Euc^5$, $S^5$, $\Hyp^5$ \\
        \rule[-9pt]{0pt}{0pt}
        $\op{U}(2)$ &
        $\Heis_5$ and $\widetilde{\op{U}(2,1)/\op{U}(2)}$ \\
        \rule[-9pt]{0pt}{0pt}
        $\op{SO}(3)_5$ &
        $\op{SL}(3,\R)/\op{SO}(3)$ and $\op{SU}(3)/\op{SO}(3)$ \\
        \rule[-18pt]{0pt}{0pt}
        $\op{SO}(2) \times \op{SO}(2)$ &
        $\SemiR{\R^4}{x-1,\,x-1,\,x+1,\,x+1}$
        and the associated bundles (Thm.~\ref{thm:full_list}(4)) \\
        \rule[-18pt]{0pt}{0pt}
        $\op{SO}(2)$ &
        $\SemiR{\R^4}{\text{2 real roots}}$ and
        $\SemiR{\R^4}{(x-1)^2,\, x+1,\, x+1}$ \\
        \rule[-12pt]{0pt}{0pt}
        $S^1_{1/2}$ &
        All line bundles over $\mathbb{F}^4$ (Thm.~\ref{thm:full_list}(5)) \\
        $S^1_1$ &
        The two unit tangent bundles (Thm.~\ref{thm:full_list}(3)), \\
        \rule[-18pt]{0pt}{0pt} &
        $\SemiR{\R^4}{x^2,\,x^2}$, and $\SemiR{(\R \times \Heis_3)}{x_3 \to x_2 \to y}$ \\
        $\{1\}$ &
        The remaining solvable Lie groups \\
    \end{tabular}\end{center}
\end{table}

The problem divides into cases by the action of
$G_p$ on the tangent space $T_p M$
(the ``linear isotropy representation'')---more specifically,
by the highest dimension of an irreducible subrepresentation $V$.

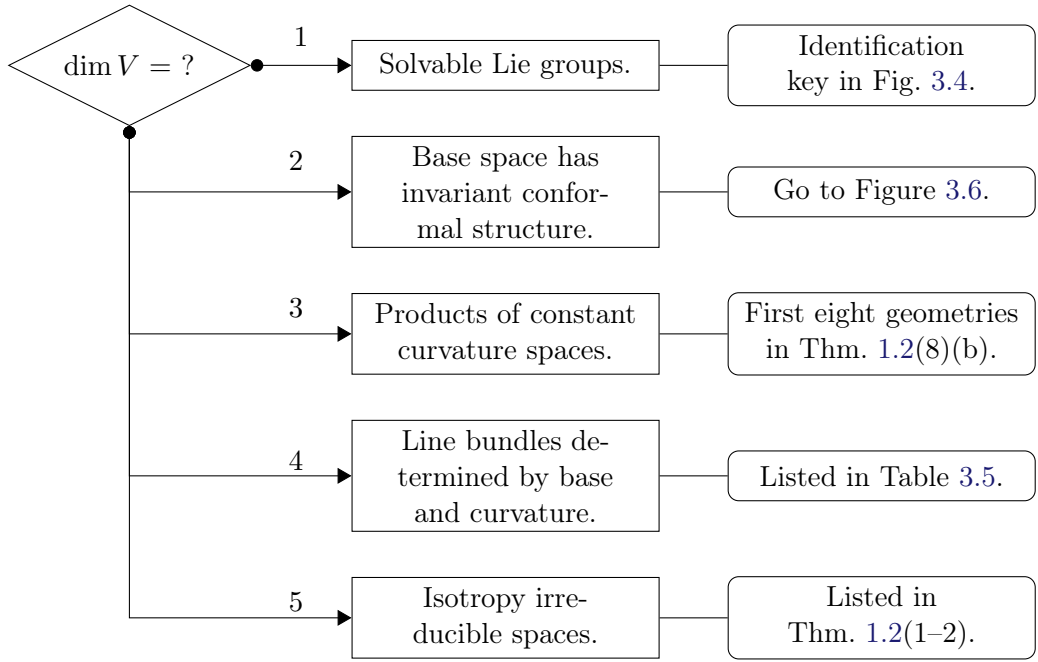
\begin{figure}[h!]
    \caption[Flowchart of the classification.]{Flowchart of the classification. Let $V$ be an irrep in $G_p \curvearrowright T_p M$ of maximal dimension.}
    \label{fig:flowchart}
\begin{center}
\begin{tikzpicture}[%
    >=triangle 60,              
    start chain=going below,    
    node distance=6mm and 50mm, 
    scale=0.8,
    ]
    \tikzset{
      base/.style={draw, on chain, on grid, align=center, minimum height=4ex},
      proc/.style={base, rectangle, text width=10em},
      test/.style={base, diamond, aspect=2, text width=5em},
      term/.style={proc, rounded corners},
      coord/.style={coordinate, on chain, on grid, node distance=6mm and 25mm},
      nmark/.style={draw, cyan, circle, font={\sffamily\bfseries}},
      it/.style={font={\small\itshape}}
    }
    \node[test] (t1) {$\dim V =$ ?};

    \node[proc, right=of t1] (v1) {Solvable Lie groups.};
    \node[proc] (v2) {Base space has invariant conformal structure.};
    \node[proc] (v3) {Products of constant curvature spaces.};
    \node[proc] (v4) {Line bundles determined by base and curvature.};
    \node[proc] (v5) {Isotropy irreducible spaces.};

    \node[term, right=of v1] (g1) {Identification key in Fig.~\ref{fig:solvable_flowchart}.};
    \node[term, right=of v2] (g2) {Go to Figure \ref{fig:fiber2_flowchart}.};
    \node[term, right=of v3] (g3) {First eight geometries in Thm.~\ref{thm:full_list}(8)(b).};
    \node[term, right=of v4] (g4) {Listed in Table \ref{table:fiber4_candidates}.};
    \node[term, right=of v5] (g5) {Listed in Thm.~\ref{thm:full_list}(1--2).};

    \path (t1.east) to node [yshift=1em]  {$1$} (v1);
        \draw[*->] (t1.east) -- (v1);
    \path (t1.south) to node [near end, yshift=0.5em] {$2$} (v2.west);
        \draw[*->] (t1.south) |- (v2);
    \path (t1.south) to node [near end, yshift=-0.9em] {$3$} (v3.west);
        \draw[*->] (t1.south) |- (v3);
    \path (t1.south) to node [near end, yshift=-2.5em] {$4$} (v4.west);
        \draw[*->] (t1.south) |- (v4);
    \path (t1.south) to node [near end, yshift=-3.7em] {$5$} (v5.west);
        \draw[*->] (t1.south) |- (v5);

    \draw (v1) -- (g1);
    \draw (v2) -- (g2);
    \draw (v3) -- (g3);
    \draw (v4) -- (g4);
    \draw (v5) -- (g5);
\end{tikzpicture}
\end{center}
\end{figure}

\begin{figure}[h!]
    \caption[Identification key for solvable geometries.
        (Duplicate of \ref{fig_thesis:solvable_flowchart})]{Identification key for solvable geometries $G = G/\{1\}$.}
    \label{fig:solvable_flowchart}
    \begin{center}\begin{tikzpicture}[%
            grow via three points={one child at (1.0,-0.7) and
            two children at (1.0,-0.7) and (1.0,-1.4)},
            growth parent anchor=west,
            edge from parent path={($(\tikzparentnode.south west)+(0.5,0)$) |- (\tikzchildnode.west)}]
        \tikzset{
            every node/.style={anchor=west},
            final/.style={draw=none},
        }
        \node {Lie algebra $\lie{g}$ is solvable}
            child { node {nilpotent}
                child { node {$4$-D abelian ideal}
                    child { node (tnil5) {$\lie{g}^4 \neq 0$}}
                    child { node (tnil4xe) {$\lie{g}^4 = 0$}}
                }
                child [missing] {}
                child [missing] {}
                child { node {no $4$-D abelian ideal}
                    child { node (tnil431) {$\lie{g}^4 \neq 0$}}
                    child { node (tnil43) {$\lie{g}^4 = 0$}}
                }
                child [missing] {}
                child [missing] {}
            }
            child [missing] {}
            child [missing] {}
            child [missing] {}
            child [missing] {}
            child [missing] {}
            child [missing] {}
            child { node {non-nilpotent}
                child { node (tr3) {nilradical $\R^3$}}
                child { node {nilradical $\R^4$}
                    child { node (tj2) {$2$ Jordan blocks}}
                    child { node (tj3) {$3$ Jordan blocks}}
                    child { node (tj4) {$4$ Jordan blocks}}
                }
                child [missing] {}
                child [missing] {}
                child [missing] {}
                child { node {nilradical $\R \oplus \lie{n}_3$}
                    child { node (tz1) {$1$-D center}}
                    child { node (tz2) {$2$-D center}}
                }
            };
        \node[right=of tj4] (geoms) {};
        \draw[->] (tr3) -- (geoms |- tr3) node {$\R^3 \rtimes \{xyz=1\}^0$};
        \draw[->] (tnil5) -- ($(geoms |- tnil5)+(0,0.25)$) node {$\SemiR{\R^4}{x^4}$};
        \draw[->] (tnil4xe) -- (geoms |- tnil4xe) node {$\SemiR{\R^4}{x^3,\, x} \cong \Nil^4 \times \Euc$};
        \draw[->] (tnil431) -- (geoms |- tnil431) node {$\SemiR{\Nil^4}{4\to3\to1}$};
        \draw[->] (tnil43) -- ($(geoms |- tnil43)+(0,-0.25)$) node {$\SemiR{\Nil^4}{3\to1}$};
        \draw[->] (tj2) -- ($(geoms |- tj2)+(0,0.25)$) node {$\SemiR{\R^4}{(x-1)^2,\, (x+1)^2}$};
        \draw[->] (tj3) -- (geoms |- tj3) node {$\SemiR{\R^4}{x^2,\, x-1,\, x+1}$};
        \draw[->] (tj4) -- ($(geoms |- tj4)+(0,-0.25)$) node {$\SemiR{\R^4}{x-a,\, x-b,\, x-c,\, x+a+b+c}$};
        \draw[->] (tz1) -- (geoms |- tz1) node {$\SemiR{(\R \times \Heis_3)}{\text{Lorentz},\, y \to x_1}$};
        \draw[->] (tz2) -- ($(geoms |- tz2)+(0,-0.25)$) node {$\Sol^4_1 \times \Euc$};
    \end{tikzpicture}\end{center}
\end{figure}
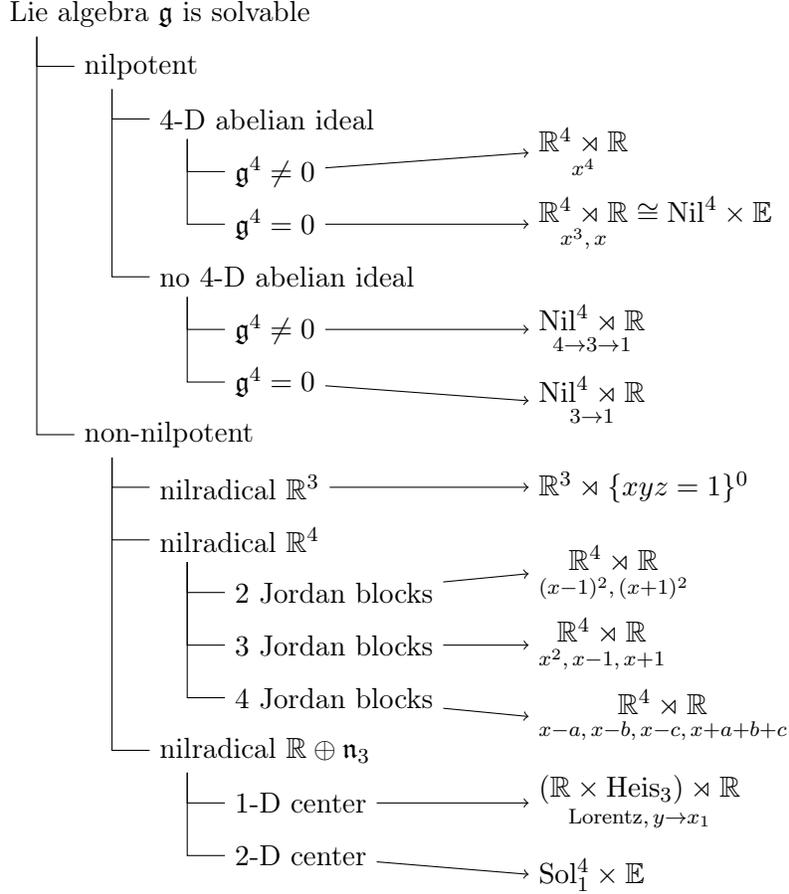

At the extremes,
one can appeal to existing classifications---the
classification of strongly isotropy irreducible homogeneous spaces by
Manturov \cite{manturov1, manturov2, manturov3, manturov1998},
Wolf \cite{wolf_irr, wolf_irr_fix}, and Kr\"amer \cite{kramer1975}
when $G_p \curvearrowright T_p M$ is irreducible ($\dim V = 5$);
and the classification of low-dimensional solvable real Lie algebras
by Mubarakzyanov \cite{muba_solvable5} and Dozias \cite{dozias}
if $G_p \curvearrowright T_p M$ is trivial ($\dim V = 1$).
These cases are handled in Part II \cite{geng2},
including the production of an identification key for trivial-isotropy
geometries (Fig.~\ref{fig:solvable_flowchart}).

Otherwise $G_p \curvearrowright T_p M$ is nontrivial and reducible
($2 \leq \dim V \leq 4$).
We classify these---the ``fibering geometries''---in Part III, starting
by proving the existence of a $G$-invariant fiber bundle structure on $M$
\cite[Prop.~\ref{iii:prop:fibering_description}]{geng3}.
The Uniformization Theorem and version of a theorem by Obata and Lelong-Ferrand
\cite[Lemma 1]{obata1973} imply the base space has an invariant conformal structure.
Beyond this common behavior,
the properties of the fibering and the relevant tools vary with the dimension of
the subrepresentation $V$, naturally suggesting the cases in
Figure \ref{fig:flowchart}.

When $\dim V = 4$, the geometries are determined by curvature and base,
in a fashion closely resembling Thurston's
treatment of $\dim V = 2$ and $\dim M = 3$
in \cite[Thm.~3.8.4(b)]{thurstonbook};
Table \ref{table:fiber4_candidates} lists the results.

\begin{table}[h!]
    \caption[Geometries with irreducible $4$-dimensional isotropy.
        (Duplicate of \ref{table_thesis:fiber4_candidates})]{Geometries with irreducible $4$-dimensional isotropy summand}
    \label{table:fiber4_candidates}
    \begin{center}\begin{tabular}{c|cc}
        Base  &
        Flat (product)  &
        Curved  \\
        \hline
        \rule[0pt]{0pt}{13pt}
        $S^4$  &
        $S^4 \times \Euc$  &
        \\
        $\Euc^4$  &
        non-maximal $\Euc^5$  &
        \\
        $\Hyp^4$  &
        $\Hyp^4 \times \Euc$  &
        \\
        $\mathbb{CP}^2$  &
        $\mathbb{CP}^2 \times \Euc$  &
        non-maximal $S^5$  \\
        $\C^2$  &
        non-maximal $\Euc^5$  &
        $\Heis_5$  \\
        $\C \Hyp^2$  &
        $\C \Hyp^2 \times \Euc$  &
        $\widetilde{\op{U}(2,1)/\op{U}(2)}$  \\
    \end{tabular}\end{center}
\end{table}

\begin{figure}[h!]
    \caption[Geometries fibering over $2$-D spaces.
        (Duplicate of \ref{fig_thesis:fiber2_flowchart})]{Classification strategy for geometries $M = G/G_p$ fibering over 2-D spaces $B$.}
    \label{fig:fiber2_flowchart}
    \begin{center}\begin{tikzpicture}[%
            >=triangle 60,              
            start chain=going below,    
            node distance=6mm and 60mm, 
            scale=0.8,
            ]
        \tikzset{
          base/.style={draw, on chain, on grid, align=center, minimum height=4ex},
          proc/.style={base, rectangle, text width=12em},
          test/.style={base, diamond, aspect=2, text width=5em},
          term/.style={proc, rounded corners},
          coord/.style={coordinate, on chain, on grid, node distance=6mm and 25mm},
          nmark/.style={draw, cyan, circle, font={\sffamily\bfseries}},
          it/.style={font={\small\itshape}}
        }
        \node[test] (tess) {$B$ has invariant metric?};
        \node[test] (tnil) {Extension of $\op{Isom}_0 \Euc^2$ by $\R^3$?};
        \node[test] (tlevi) {Levi action nontrivial?};
        \node[term] (gprod) {Products and associated bundles.};

        \node[term, right=of tess] (gess)
            {$T^1 \Hyp^3$ and some non-nilpotent solvable Lie groups.};
        \node[term, right=of tnil] (gnil) {Some nilpotent Lie groups.};
        \node[term, right=of tlevi] (glevi) {$T^1 \Euc^{1,2}$ and the line bundles over $\mathbb{F}^4$.};

        \path (tess.east) to node [near start, yshift=1em] {no} (gess);
            \draw[*->] (tess.east) -- (gess);
        \path (tnil.east) to node [near start, yshift=1em] {yes} (gnil);
            \draw[*->] (tnil.east) -- (gnil);
        \path (tlevi.east) to node [near start, yshift=1em] {yes} (glevi);
            \draw[*->] (tlevi.east) -- (glevi);

        \path (tess.south) to node [anchor=west, xshift=0.5em] {yes---$\tilde{G}$ is an extension of $\op{Isom}_0 B$} (tnil);
            \draw[o->] (tess.south) -- (tnil);
        \path (tnil.south) to node [anchor=west, xshift=0.5em] {no---then $\tilde{G}$ is a split extension} (tlevi);
            \draw[o->] (tnil.south) -- (tlevi);
        \path (tlevi.south) to node [anchor=west, xshift=0.5em] {no---$\tilde{G}$ is a direct product} (gprod);
            \draw[o->] (tlevi.south) -- (gprod);
    \end{tikzpicture}\end{center}
\end{figure}
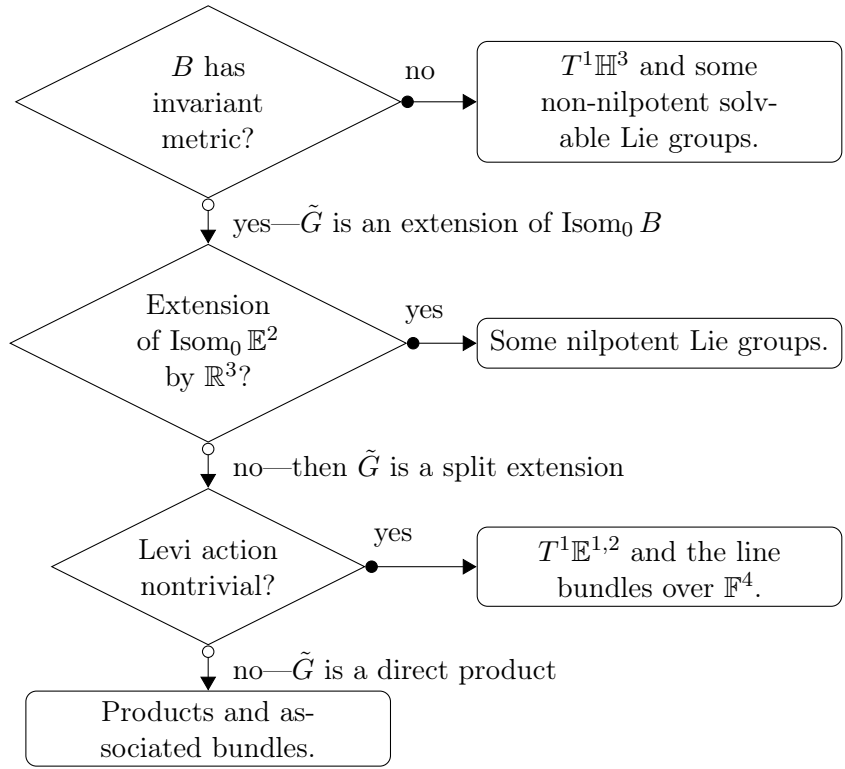

Otherwise, we work systematically with $G$-invariant fiber bundle structures
by recasting the problem as an extension problem for the Lie algebra of $G$
and solving it with the help of Lie algebra cohomology.
Over $3$-dimensional base spaces there happen to be only products;
but over $2$-dimensional base spaces a daunting array of possibilities
requires some attempt to organize the problem, summarized in
Figure \ref{fig:fiber2_flowchart}.

\section{Related work} \label{sec:whatelse}

\paragraph{Classification of compact homogeneous spaces.}

Gorbatsevich has produced classification results for
\emph{compact} homogeneous spaces $M$
by using a fiber bundle described in \cite[\S II.5.3.2]{onishchik1}
whose fibers have compact transformation group, whose base is aspherical,
and whose total space is a finite cover of $M$.
The classification is complete in dimension
up to $5$ in general, in dimension $6$ up to finite covers, and in dimension $7$
in the aspherical case \cite{gorbatsevich2012}.

Another approach would be to group the problem by the number of isotropy
summands. The \emph{Riemannian} homogeneous spaces with irreducible isotropy
were classified by Manturov \cite{manturov1, manturov2, manturov3, manturov1998},
Wolf \cite{wolf_irr, wolf_irr_fix}, and Kr\"amer \cite{kramer1975};
and the \emph{compact} \emph{Riemannian} homogeneous spaces with two isotropy summands
were classified by Dickinson and Kerr in \cite{dickinsonkerr}.

\paragraph{Classification of naturally reductive spaces.}

The \keyword{naturally reductive} Riemannian homogeneous spaces $G/G_p$---those
whose geodesics through $p$ are orbits of $1$-parameter subgroups tangent to
the representation complementary to $\tanalg G_p$ in $G_p \curvearrowright \tanalg G$
(see e.g.\ \cite[\S{X.3}]{kobayashinomizu2})---have been
classified in dimension $6$ by Agricola, Ferreira, and Friedrich
\cite{aff2015}; and in lower dimensions by work of Kowalski and Vanhecke
(see \cite[\S{}6]{kowalski1996} for a summary).

The case of dimension $5$, in \cite[Thm.~2.1]{kowalski1985},\footnote{
    The changes in the corrected version \cite[6.4]{kowalski1996}
    appear to amount to (1) changing ``symmetric'' and ``decomposable''
    to ``locally symmetric'' and ``locally decomposable'' and (2)
    changing the rational parameter to a real parameter in the Type II
    family (Heisenberg-group-like bundles).
}
shares the following features with the classification of geometries.
\begin{enumerate}
    \item Everything with $\op{SU}(2)$ isotropy
        is realized by a homogeneous space with $\op{U}(2)$ isotropy
        \cite[main result (b)]{kowalski1991}.
    \item The associated bundle geometries appear
        as indecomposable naturally reductive spaces.
\end{enumerate}
The differences between geometries and naturally reductive spaces bear
mentioning as well:
\begin{enumerate}
    \item Naturally reductive spaces need not be maximal as geometries,
        as demonstrated by non-maximal realizations of $S^3 \cong S^3 \rtimes S^1 / S^1$
        and $S^5 \cong \op{SU}(3)/\op{SU}(2)$.
    \item Some geometries---particularly those with trivial isotropy---may not
        be realizable by naturally reductive spaces.
        In $3$ dimensions, there is just $\Sol^3$;
        and in $4$ dimensions, there are $\mathbb{F}^4$ and the four
        solvable Lie group geometries other than
        $\Heis_3 \times \Euc$ and $\Euc^4$.
        In $5$ dimensions, there are
        the unit tangent bundles $T^1 \Hyp^3$ and $T^1 \Euc^{1,2}$,
        the line bundles over $\mathbb{F}^4$,
        the products involving $\Sol^3$,
        and any solvable Lie group geometries
        that are not $\Euc^5$, $\Heis_5$,
        or a product involving $\Heis_3$.
\end{enumerate}
A chart in \cite[5.1]{kowalski1996} summarizes the relations between
several other classes of spaces.

\paragraph{Other geometric structures.}
One can replace the assumption of an invariant Riemannian metric
(compact point stabilizers) with other structures.
A number of difficulties may result from this:
geodesic completeness may no longer coincide with other notions
of completeness (e.g.\ in \cite[Thm.~2.1]{dz2010});
isotropy representations may fail to be semisimple or faithful;
and point stabilizers may fail to act faithfully on tangent spaces,
necessitating techniques like those of \cite{capmelnick}.
In spite of these challenges, some results are known, such as:
\begin{itemize}
    \item Using conformal structures without relaxing other assumptions
        yields only $S^n$ and $\Euc^n$:
        a manifold whose conformal automorphism group acts transitively
        with the identity component preserving no Riemannian metric
        is conformally equivalent to one of the two, by theorems of
        Obata \cite[Lemma 1]{obata1973} and Lafontaine \cite[Thm.~D.1]{lafontaine}.
    \item The interaction of complex structures and $4$-dimensional geometries
        was investigated by Wall in \cite{wall86}; and almost-complex structures
        on homogeneous spaces up to dimension $6$ are classified by
        Alekseevsky, Kruglikov, and Winther in
        \cite{winther} with the additional assumption that point stabilizers
        are semisimple.
    \item Loosely analogous to almost-complex structures
        are $5$-manifolds whose structure group can be reduced to $\op{SO}(3)_5$.
        The classification of such structures
        in the integrable case by Bobie\'nski and Nurowski
        \cite[Thm.~4.7]{bn2005} offers an alternative path to the classification
        of isotropy-irreducible geometries,
        and there are some further classification results in the non-integrable case
        by Chiossi and Fino in \cite{fc2007} and
        by Agricola, Becker-Bender, and Friedrich in \cite{abf2011}.
    \item The pseudo-Riemannian geometries were
        classified in dimension $3$ by Dumitrescu and Zeghib in \cite{dz2010};
        and the pseudo-Riemannian naturally reductive spaces were classified in
        dimension $4$ by Batat, L\'opez, and Mar\'ia \cite{blm2015}.
\end{itemize}

\bibliographystyle{amsalpha}
\bibliography{main}

\providecommand{\bysame}{\leavevmode\hbox to3em{\hrulefill}\thinspace}
\providecommand{\MR}{\relax\ifhmode\unskip\space\fi MR }
\providecommand{\MRhref}[2]{%
  \href{http://www.ams.org/mathscinet-getitem?mr=#1}{#2}
}
\providecommand{\href}[2]{#2}
\begin{thebibliography}{PSWZ76}

\bibitem[ABBF11]{abf2011}
I.~Agricola, J.~Becker-Bender, and T.~Friedrich, \emph{On the topology and the
  geometry of so(3)-manifolds}, Annals of Global Analysis and Geometry
  \textbf{40} (2011), 67--84, arXiv: 1010.0260.

\bibitem[AFF15]{aff2015}
I.~Agricola, A.~C. Ferreira, and T.~Friedrich, \emph{The classification of
  naturally reductive homogeneous spaces in dimensions n $\leq$ 6},
  Differential Geometry and its Applications \textbf{39} (2015), 59--92.

\bibitem[AKW14]{winther}
D.V. Alekseevsky, B.~Kruglikov, and H.~Winther, \emph{Homogeneous almost
  complex structures in dimension 6 with semi-simple isotropy}, Annals of
  Global Analysis and Geometry \textbf{46} (2014), 361--387, arXiv: 1401.8187.

\bibitem[Bar65]{barden}
D.~Barden, \emph{Simply connected five-manifolds}, Annals of Mathematics
  \textbf{82} (1965), no.~3, 365--385.

\bibitem[BFNT13]{boza}
L.~Boza, E.~M. Fedriani, J.~Nunez, and A.~F. Tenorio, \emph{A historical review
  of the classifications of {Lie} algebras}, Rev. Un. Mat. Argentina
  \textbf{54} (2013), no.~2.

\bibitem[BG02]{boyergalicki}
C.~P. Boyer and K.~Galicki, \emph{Rational homology 5-spheres with positive
  {R}icci curvature}, Mathematical Research Letters \textbf{9} (2002), no.~4,
  521--528, arXiv: math/0203048.

\bibitem[BLM15]{blm2015}
W.~Batat, M.~C. López, and E.~R. María, \emph{Four-dimensional naturally
  reductive pseudo-{Riemannian} spaces}, Differential Geometry and its
  Applications \textbf{41} (2015), 48--64.

\bibitem[BN07]{bn2005}
M.~Bobieński and P.~Nurowski, \emph{Irreducible so(3)-geometries in dimension
  five}, Journal für die reine und angewandte Mathematik \textbf{605} (2007),
  51--93, arXiv: math/0507152.

\bibitem[DK08]{dickinsonkerr}
W.~Dickinson and M.~M. Kerr, \emph{The geometry of compact homogeneous spaces
  with two isotropy summands}, Annals of Global Analysis and Geometry
  \textbf{34} (2008), no.~4, 329--350.

\bibitem[Doz63]{dozias}
J.~Dozias, \emph{Sur les alg\`ebres de {L}ie r\'esolubles r\'eelles de
  dimension inf\'erieure ou \'egale \`a 5}, Ph.D. thesis, Facult\'e des
  Sciences de Paris, 1963.

\bibitem[Dyn00a]{dynkin2000maximal}
E.~B. Dynkin, \emph{Maximal subgroups of the classical groups}, Selected papers
  of E. B. Dynkin with commentary, Amer. Math. Soc., Providence, RI, 2000,
  pp.~37--170.

\bibitem[Dyn00b]{dynkin2000semisimple}
\bysame, \emph{Semisimple subalgebras of semisimple {L}ie algebras}, Selected
  papers of E. B. Dynkin with commentary, Amer. Math. Soc., Providence, RI,
  2000, pp.~175--308.

\bibitem[DZ10]{dz2010}
S.~Dumitrescu and A.~Zeghib, \emph{Géométries {Lorentziennes} de dimension 3:
  classification et complétude}, Geometriae Dedicata \textbf{149} (2010),
  no.~1, 243--273.

\bibitem[FC07]{fc2007}
A.~Fino and S.~Chiossi, \emph{Nearly integrable so(3) structures on
  5-dimensional lie groups}, J. Lie Theory \textbf{17} (2007), no.~3, 539--562.

\bibitem[Fil83]{filipk}
R.~Filipkiewicz, \emph{Four dimensional geometries}, Ph.D. thesis, University
  of Warwick, 1983.

\bibitem[Gen16a]{geng2}
A.~Geng, \emph{5-dimensional geometries ii: the non-fibered geometries}, In
  preparation, 2016.

\bibitem[Gen16b]{geng3}
\bysame, \emph{5-dimensional geometries iii: the fibered geometries}, In
  preparation, 2016.

\bibitem[Gor77]{gorbatsevich1977}
V.~V. Gorbatsevich, \emph{Three-dimensional homogeneous spaces}, Siberian
  Mathematical Journal \textbf{18} (1977), no.~2, 200--210 (en).

\bibitem[Gor12]{gorbatsevich2012}
\bysame, \emph{On compact aspherical homogeneous manifolds of dimension
  $\leq$7}, Mathematical Notes \textbf{92} (2012), no.~1-2, 186--196 (en).

\bibitem[GOV93]{onishchik1}
V.~V. Gorbatsevich, A.~L. Onishchik, and E.~B. Vinberg, \emph{Lie groups and
  {L}ie algebras {I}: Foundations of {L}ie theory, {L}ie transformation
  groups}, Springer-Verlag, 1993.

\bibitem[Hil02]{hillman}
J.~A. Hillman, \emph{Four-manifolds, geometries and knots}, University of
  Warwick, Mathematics Institute, 2002.

\bibitem[Hus94]{husemoller}
D.~Husemöller, \emph{Fibre bundles}, 3rd ed., Graduate texts in mathematics,
  no.~20, Springer-Verlag, New York, 1994.

\bibitem[HZ96]{heintze1996}
E.~Heintze and W.~Ziller, \emph{Isotropy irreducible spaces and
  s-representations}, Differential Geometry and its Applications \textbf{6}
  (1996), no.~2, 181--188.

\bibitem[Ish55]{ishihara1955}
Shigeru Ishihara, \emph{Homogeneous {R}iemannian spaces of four dimensions.},
  Journal of the Mathematical Society of Japan \textbf{7} (1955), no.~4,
  345--370 (EN). \MR{82717}

\bibitem[KN69]{kobayashinomizu2}
S.~Kobayashi and K.~Nomizu, \emph{Foundations of differential geometry},
  Interscience Tracts in Pure and Applied Mathematics, vol.~2, Interscience
  Publishers (a division of John Wiley \& Sons), New York-London, 1969.

\bibitem[KPV96]{kowalski1996}
O.~Kowalski, F.~Pr\"ufer, and L.~Vanhecke, \emph{D’{Atri} spaces}, Topics in
  geometry, Springer, 1996, pp.~241--284.

\bibitem[Kr{\"a}75]{kramer1975}
M.~Kr{\"a}mer, \emph{Eine klassifikation bestimmter untergruppen kompakter
  zusammenhängender {L}iegruppen}, Communications in Algebra \textbf{3}
  (1975), no.~8, 691--737.

\bibitem[KV85]{kowalski1985}
O.~Kowalski and L.~Vanhecke, \emph{Classification of five-dimensional naturally
  reductive spaces}, Math. {Proc}. {Camb}. {Phil}. {Soc}, vol.~97, Cambridge
  Univ Press, 1985, pp.~445--463.

\bibitem[KV91]{kowalski1991}
O.~Kowalski and L.~Vanhecke, \emph{Riemannian manifolds with homogeneous
  geodesics}, Bollettino della Unione Matemàtica Italiana. Serie VII. B
  \textbf{5} (1991), no.~1, 189--246.

\bibitem[Laf88]{lafontaine}
J.~Lafontaine, \emph{The theorem of {L}elong-{F}errand and {O}bata}, Conformal
  geometry, Springer, 1988, pp.~93--103.

\bibitem[Man61a]{manturov1}
O.~V. Manturov, \emph{Homogeneous asymmetric {R}iemannian spaces with an
  irreducible group of rotations}, Dokl. Akad. Nauk. SSSR \textbf{141} (1961),
  no.~4, 792--795.

\bibitem[Man61b]{manturov2}
\bysame, \emph{Riemannian spaces with orthogonal and symplectic groups of
  motions and an irreducible group of motions}, Dokl. Akad. Nauk. SSSR
  \textbf{141} (1961), no.~5, 1034--1037.

\bibitem[Man66]{manturov3}
\bysame, \emph{Homogeneous {R}iemannian manifolds with irreducible isotropy
  group}, Trudy Sem. Vector and Tensor Analysis \textbf{13} (1966), 68--145.

\bibitem[Man98]{manturov1998}
\bysame, \emph{Homogeneous {R}iemannian spaces with irreducible rotation
  group}, Tensor and Vector Analysis. Geometry, Mechanics and Physics, Gordon
  and Breach, Amsterdam (1998), 101--192.

\bibitem[Mat02]{matsumoto_morse}
Y.~Matsumoto, \emph{An introduction to morse theory}, vol. 208, American
  Mathematical Soc., 2002.

\bibitem[Mos50]{mostow1950}
G.~D. Mostow, \emph{The extensibility of local {Lie} groups of transformations
  and groups on surfaces}, Annals of Mathematics \textbf{52} (1950), no.~3,
  606--636.

\bibitem[Mos55]{mostow1955}
\bysame, \emph{On covariant fiberings of {K}lein spaces}, American Journal of
  Mathematics \textbf{77} (1955), no.~2, 247--278.

\bibitem[Mos62]{mostow_covariant_1962}
\bysame, \emph{Covariant fiberings of {Klein} spaces, ii}, American Journal of
  Mathematics \textbf{84} (1962), no.~3, 466--474.

\bibitem[Mub63]{muba_solvable5}
G.~M. Mubarakzyanov, \emph{Classification of real structures of {L}ie algebras
  of fifth order}, Izv. Vyssh. Uchebn. Zaved. Mat. (1963), no.~3, 99--106.
  \MR{155871}

\bibitem[Oba73]{obata1973}
M.~Obata, \emph{A non-compact {Riemannian} manifold admitting a transitive
  group of conformorphisms}, Tohoku Mathematical Journal \textbf{25} (1973),
  no.~4, 553--556 (EN). \MR{0346704}

\bibitem[Ott09]{ottenburger2009}
S.~Ottenburger, \emph{A diffeomorphism classification of 5-and 7-dimensional
  non-simply-connected homogeneous spaces}, Ph.D. thesis, Bonn, Univ., Diss.,
  2009, 2009.

\bibitem[OV00]{onishchik2}
A.~L. Onishchik and E.~B. Vinberg, \emph{{L}ie groups and {L}ie algebras {II}:
  Discrete subgroups of {L}ie groups and cohomologies of {L}ie groups and {L}ie
  algebras}, Encyclopaedia of Mathematical Sciences, vol.~21, Springer-Verlag,
  2000.

\bibitem[PSWZ76]{patera}
J.~Patera, R.~T. Sharp, P.~Winternitz, and H.~Zassenhaus, \emph{Invariants of
  real low dimension {Lie} algebras}, Journal of Mathematical Physics
  \textbf{17} (1976), no.~6, 986--994.

\bibitem[Sco83]{scott}
P.~Scott, \emph{The geometries of 3-manifolds}, Bulletin of the London
  Mathematical Society \textbf{15} (1983), no.~5, 401--487. \MR{705527
  (84m:57009)}

\bibitem[Sha00]{sharpe}
R.~W. Sharpe, \emph{Differential geometry: {Cartan}'s generalization of
  {Klein}'s {Erlangen} program}, Springer New York, December 2000 (en).

\bibitem[{\v{S}}W12]{snobl2012}
L.~{\v{S}}nobl and P.~Winternitz, \emph{Solvable {L}ie algebras with {B}orel
  nilradicals}, Journal of Physics A: Mathematical and Theoretical \textbf{45}
  (2012), no.~9, 095202.

\bibitem[Thu97]{thurstonbook}
W.~P. Thurston, \emph{Three-dimensional geometry and topology}, vol.~1,
  Princeton University Press, 1997.

\bibitem[vM13]{capmelnick}
A.~\v{C}ap and K.~Melnick, \emph{Essential {Killing} fields of parabolic
  geometries}, Indiana University Mathematics Journal \textbf{62} (2013),
  no.~6, 1917--1953, arXiv: 1208.5510.

\bibitem[Wal86]{wall86}
C.~T.~C. Wall, \emph{Geometric structures on compact complex analytic
  surfaces}, Topology \textbf{25} (1986), no.~2, 119--153.

\bibitem[Wol68]{wolf_irr}
J.~A. Wolf, \emph{The geometry and structure of isotropy irreducible
  homogeneous spaces}, Acta Mathematica \textbf{120} (1968), no.~1, 59--148
  (en).

\bibitem[Wol84]{wolf_irr_fix}
\bysame, \emph{Erratum to: {The} geometry and structure of isotropy irreducible
  homogeneous spaces}, Acta Mathematica \textbf{152} (1984), no.~1, 141--142.

\end{thebibliography}

\end{document}